\documentclass[12pt,reqno]{article}
\usepackage{amssymb}
\usepackage{amsthm}
\usepackage{amsmath}
\newtheorem{theorem}{Theorem}[section]

\newtheorem{example}[theorem]{Example}

\def\be{\begin{equation}}
\def\ee{\end{equation}}
\def\bea{\begin{eqnarray}}
\def\eea{\end{eqnarray}}
\begin{document}
\begin{center} \Large{\bf Invariant characterization of third-order ODEs $u'''=f(x,u,u',u'')$ that admit a five-dimensional point symmetry Lie algebra}
\end{center}
\medskip
\begin{center}
Ahmad Y. Al-Dweik$^*$, M. T. Mustafa$^{**}$ and
F. M. Mahomed$^{***,}$\\

{$^*$Department of Mathematics \& Statistics, King Fahd University
of Petroleum and Minerals, Dhahran 31261, Saudi Arabia}\\
{$^{**}$Department of Mathematics, Statistics and Physics, Qatar
University, Doha, 2713, State of Qatar}\\
{$^{***}$School of Computer Science and Applied Mathematics, DST-NRF
Centre of Excellence in Mathematical and Statistical Sciences,
University of the Witwatersrand, Johannesburg, Wits 2050,  South Africa\\
}

aydweik@kfupm.edu.sa, tahir.mustafa@qu.edu.qa and
Fazal.Mahomed@wits.ac.za
\end{center}
\begin{abstract}
The Cartan equivalence method is applied to provide an invariant
characterization of the third-order ordinary differential equation
$u'''=f(x,u,u',u'')$ which admits a five-dimensional point
symmetry Lie algebra. The invariant characterization is given in
terms of the function $f$ in a compact form. A simple procedure to
construct the equivalent canonical form by use of an obtained
constant invariant is also presented. We also show how one obtains
the point transformation that does the reduction to linear form.
Moreover, some applications are provided.
\end{abstract}
\newpage
\section{Introduction}
Both practical and algebraic linearization criteria for scalar
second-order ordinary differential equations (ODEs) by means of
invertible point transformations were first obtained by Lie
(\cite{lie, lie1}). Lie \cite{lie1} showed that the most general
form of a scalar second-order ODE which is reducible to a linear
equation via maps of the independent and dependent variables is
at most cubic in the first-order derivative and he moreover
provided a practical linearization test in terms of the
coefficients of the cubic equation. Lie \cite{lie} also worked out
the algebraic criteria for linearizability for such equations. Any
scalar second-order ODE with the maximum number which is eight of point
symmetries is linearizable. Lie \cite{lie} further showed that if
a second-order ODE admits a two-dimensional Lie algebra of point
symmetries which is of rank 1, then it is linearizable by point transformation (see e.g.
Mahomed \cite{mah1}).

Tress\'e proved that the scalar second-order ODE \be
y''=f(x,y,y')\label{e1} \ee is linearizable by means of a point
transformation if and only if the Tress\'e \cite{tre} relative
invariants
\begin{eqnarray}
I&=&f_{y'y'y'y'},\nonumber\\
J&=&\frac{d^2}{dx^2} f_{y'y'}-4\frac{d}{dx} f_{y'y}-3f_yf_{y'y'}+
6f_{yy}+f_{y'}(4f_{y'y}-\frac{d}{dx} f_{y'y'})\label{e2}
\end{eqnarray}
both vanish identically for the said equation (\ref{e1}). Thus Tress\'e
\cite{tre} considered the linearization problem for scalar
second-order ODEs in terms of the Tress\'e relative invariants of
the equivalence group of point transformations. The setting to
zero of the Tress\'e invariants is equivalent to the compatibility
of the over-determined Lie  conditions for linearization
\cite{mah2}.

Another method for studying the equivalence problem, in particular
the linearzation problem, for ODEs was considered by Cartan
\cite{car} called the Cartan equivalent method (see
\cite{Gardner1989, Olver1995}). This approach associates a
geometric structure with each differential equation. Grissom et.
al. \cite{gri} used the Cartan equivalence method to arrive at the
Lie invariant criteria for linearization for scalar second-order
ODEs.

Yet another approach is a geometric one, viz. that of projection of the
geodesic equations by dimension one as proposed in \cite{mah3}.
This enables a geometric re-derivation of Lie's linearization
conditions for a scalar second-order ODE. Furthermore, it is shown
how the point transformations for reduction to a linear equation
can be constructed in the higher space and by utilization of the
coefficients of the original ODE.

Our focus here is to study the linearization problem via invertible maps for scalar
third-order ODEs which admit a five dimensional point symmetry
algebra using the Cartan approach. Thus we firstly review works relevant to scalar third-order ODEs.
Mahomed and Leach \cite{mah4} found the algebraic criteria for
linearization for scalar $n$th-order ($n>2$) ODEs. For scalar
linear third-order ODEs they deduced three forms. The
Laguerre-Forsyth (see \cite{mah4}) canonical form for such
third-order equations is given by \be u'''+au=0\label{e3} \ee
where $a=a(x)$. If $a=0$, (\ref{e3}) has the maximal symmetry
Lie algebra of dimension 7. Otherwise the equation (\ref{e3}) admits a
five- or four-dimensional symmetry algebra.

After the pioneering works of Lie \cite{lie,lie1} and Tress\'e
\cite{tre}, there has been a renewal of interest in invariant
linearization criteria for ODEs. Here we pay attention to scalar
third-order ODEs. Chern \cite{che} provided a major impetus in the
solution of the linearization problem of scalar third-order
equations by means of contact transformations by invocation of the Cartan
equivalence method. He derived conditions for equivalence to (\ref{e3}) for
$a=0$ and $a=1$. The linearization conditions are in terms of
geometric invariants of contact transformations. Neut and Petitot
\cite{neu} investigated conditions on equivalence to (\ref{e3}).
We mention and review these below. Grebot \cite{gre} also focused on the
linearization of third-order ODEs. However, this was via fibre
preserving or a restricted class of point transformations.
Ibragimov and Meleshko \cite{ibr} studied the linearization
problem for scalar third-order ODEs by means of point and contact
transformations. The conditions on the linearizing transformations
are provided in their works as well.  They invoke the Laguerre-Forsyth canonical form
for reduction. Conditional invariant linearization criteria for
scalar third-order ODEs are given in Mahomed and Qadir
\cite{mah5}. These conditions are for third-order ODEs which are
solvable in terms of two arbitrary constants.

Our main purpose in this work is to study the linearization
problem via invertible transformations for scalar third-order ODE $u'''=f(x,u,u',u'')$ by the
Cartan equivalence method which enables reduction to the canonical form with
five symmetries
and to provide compact criteria in terms of $f$. The case of seven point
symmetries ($k=0$) was addressed by Al-Dweik \cite{Dweik2016} (see
below).

Neut and Petitot \cite{neu} proved that the necessary and
sufficient conditions for linearization of the third-order ODE $
u''' = f(x,u,u',u'')$ to the normal form $u''' = 0$  under {\it
contact transformation} are vanishing of the
W$\ddot{\textrm{u}}$nschmann relative invariants (\ref{wer2}) as
stated in the next theorem.

\begin{theorem}\cite{neu}
The equation $ u''' = f(x,u,u',u'')$ is equivalent to the normal
form $u''' = 0$ with seven point symmetries under {\rm contact
transformations} if and only if the relative invariants
\begin{equation}\label{wer2}
\begin{array}{ll}
I_1=f_{u'',u'',u'',u''}\\
I_2 =4\,f_{u''}^3  + 18\,f_{u''} \left( {f_{u'}  - \dot{D}_x f_{u''}} \right) + 9\,\dot{D}_x^2 f_{u''}  + 54\,f_u  - 27\,\dot{D}_x f_{u'}\\
\end{array}
\end{equation}
both vanish identically, where $I_2$ is the well-known
W$\ddot{\textrm{u}}$nschmann relative invariant \cite{neu,
Wünschmann}.
\end{theorem}
Invariant characterization of third-order ODEs $u'''=f(x,u,u',u'')$ which admit a seven-dimensional
point symmetry algebra was given in terms of the function $f$ in a
compact form in the following theorem. In the sequel, we denote
$u', u''$ by $p, q$, respectively.
\begin{theorem}\cite{Dweik2016}
The necessary and sufficient conditions for equivalence of a
third-order equation $ u''' = f(x,u,u',u'')$ to the canonical
form $u''' = 0$ with seven symmetries via {\rm point
transformation} are the identically vanishing of the system of
relative invariants
\begin{equation}\label{c1}
\begin{array}{ll}
{f_{q,q,q} }  \\
{f_{q,q} ^2  + 6\,f_{p,q,q} }  \\
4\,f_q ^3  + 18\,f_q \left( {f_p  - \dot{D}_x f_q} \right) + 9\,\dot{D}_x^2 f_q  + 54\,f_u  - 27\,\dot{D}_x f_p\\
f_{q,q} \left( {f_q ^2+9\,f_p  - 3\,\dot{D}_x f_q}\right) - 9\,f_{p,p}  + 18\,f_{u,q}  - 6\,f_q f_{p,q},\\
\end{array}
\end{equation}
where $ \dot{D}_x = \frac{\partial } {{\partial x}} +
p\frac{\partial }{{\partial u}}+q\frac{\partial} {{\partial
p}}+f\frac{\partial} {{\partial q}}$.
\end{theorem}
Our aim in this paper is to give the necessary and sufficient
conditions  for equivalence of third-order equations $u''' =
f(x,u,u',u'')$ via {\it point transformations} to the canonical form with five symmetries, in
terms of $f$.

Mahomed and Leach \cite{mah4} showed that a scalar third-order ODE
with a 5-dimensional symmetry algebra is linearizable via a point
transformation and is equivalent to the linear form
\begin{equation}\label{a1}
\begin{array}{cc}
u'''=ku'+lu,& k,l(\ne0)=\,{\rm constant}.
\end{array}
\end{equation}
The transformation $x=\frac{t}{l^{\frac{1}{3}}}$ maps the
canonical from (\ref{a1}) to
\begin{equation}\label{a3}
\begin{array}{cc}
\frac{d^3}{dt^3}u=\frac{k}{l^{\frac{2}{3}}}\frac{d}{dt}u+u.
\end{array}
\end{equation}
Therefore, the canonical form for third-order equations with a
5-dimensional point symmetry algebra can be simplified further to the
following canonical form
\begin{equation}\label{a4}
\begin{array}{ll}
u'''=s~u'+u,&s=\,{\rm constant},\\
\end{array}
\end{equation}
In this paper, we consider the canonical form (\ref{a4}) with five
point symmetries instead of the Laguerre-Forsyth canonical form
for third-order ODEs. The reason is that the canonical form
(\ref{a4}) has constant coefficients, while the Laguerre-Forsyth
form may have variable coefficient $a(x)$ for third-order ODEs
with five point symmetries. For example $u'''=\frac{1}{x^6}u$ has
the following five point symmetries
\begin{equation}
\begin{array}{lll}
X_1=u\frac{\partial}{\partial u},& X_2=x^2\frac{\partial}{\partial x}+2xu\frac{\partial}{\partial u},& X_3=x^2e^{-\frac{1}{x}}\frac{\partial}{\partial u},\\
X_4=x^2 e^{\frac{1}{2x}}\cos
\left({\frac{\sqrt{3}}{2x}}\right)\frac{\partial}{\partial u},&X_5=x^2 e^{\frac{1}{2x}}\sin\left({\frac{\sqrt{3}}{2x}}\right)\frac{\partial}{\partial u}.& \\
\end{array}
\end{equation}
Also it is important to mention here that apart from not utilizing
the Laguerre-Forsyth canonical form as in \cite{ibr}, we also wish
to distinguish the five symmetry linear canonical form and to
provide compact criteria in terms of $f(x,u,u',u'')$ as well as to
utilize the Cartan equivalence method whereas the work \cite{ibr}
used the direct method. Moreover, by use of the Cartan method, we
for the first time provide how one obtains the point
transformation to the five symmetry linear canonical form via the
Cartan approach.
\section{Application of Cartan's equivalence method for third-order ODEs with five point symmetries}
Let $x:=(x,u,p=u',q=u'')\in \mathbb{R}^4$ be local coordinates
of $J^2$, the space of the second order jets. In local
coordinates, the equivalence of
\begin{equation}\label{b0}
\begin{array}{cc}
u'''=f(x,u,u',u''), & \bar{u}'''=\bar{f}(\bar{x},\bar{u},\bar{u}',\bar{u}''),\\
\end{array}
\end{equation}
under a point transformation
\begin{equation}\label{ccc}
\bar{x}=\phi \left( x,u \right),~\bar{u} =\psi \left( x,u  \right),~~~\phi_x\psi_u   -  \phi_u \psi_x\neq0,\\
\end{equation}
is expressed as the local equivalence problem for the $G$-structure
\begin{equation}\label{b2}
\Phi^*\left(%
\begin{array}{c}
  \bar{\omega}^1 \\
  \bar{\omega}^2 \\
  \bar{\omega}^3 \\
  \bar{\omega}^4 \\
\end{array}%
\right)=\left(%
\begin{array}{cccc}
  a_1 & 0 & 0 & 0 \\
  a_2 & a_3 & 0 & 0 \\
  a_4 & a_5 & a_6 & 0 \\
  a_7 & 0 & 0 & a_8 \\
\end{array}%
\right) \left(%
\begin{array}{c}
  \omega^1 \\
  \omega^2 \\
  \omega^3 \\
  \omega^4 \\
\end{array}%
\right),
\end{equation}
where
\begin{equation}\label{b1}
\begin{array}{llll}
\bar{\omega}^1=d\bar{u}-\bar{p} d\bar{x}, & \bar{\omega}^2=d\bar{p}-\bar{q} d \bar{x}, & \bar{\omega}^3=d\bar{q}-\bar{f} d \bar{x}, & \bar{\omega}^4= d \bar{x},\\
\omega^1=du-p d x, & \omega^2=dp-q d x, & \omega^3=dq-f d x, & \omega^4= d x.\\
\end{array}
\end{equation}
One can evaluate the functions $a_i=a_i(x,u,p,q), i=1..8,$
explicitly. For instance, $a_1=\frac{\phi_x\psi_u   - \phi_u
\psi_x}{D_x \phi}.$

Now, let us define $\theta$ to be the lifted coframe with an
eight-dimensional group $G$
\begin{equation}\label{b21}
\left(%
\begin{array}{c}
  \theta^1 \\
  \theta^2 \\
  \theta^3 \\
  \theta^4 \\
\end{array}%
\right)=\left(%
\begin{array}{cccc}
  a_1 & 0 & 0 & 0 \\
  a_2 & a_3 & 0 & 0 \\
  a_4 & a_5 & a_6 & 0 \\
  a_7 & 0 & 0 & a_8 \\
\end{array}%
\right) \left(%
\begin{array}{c}
  \omega^1 \\
  \omega^2 \\
  \omega^3 \\
  \omega^4 \\
\end{array}%
\right).
\end{equation}
Cartan's method, when applied to this equivalence problem, leads
to an ${e}$-structure, which is invariantly associated to the
given equation.

The first structure equation is
\begin{equation}\label{b3}
d\left(%
\begin{array}{c}
  \theta^1 \\
  \theta^2 \\
  \theta^3 \\
  \theta^4 \\
\end{array}%
\right)=\left(%
\begin{array}{cccc}
  \alpha_1 & 0 & 0 & 0 \\
  \alpha_2 & \alpha_3 & 0 & 0 \\
  \alpha_4 & \alpha_5 & \alpha_6 & 0 \\
  \alpha_7 & 0 & 0 & \alpha_8 \\
\end{array}%
\right)\wedge \left(%
\begin{array}{c}
  \theta^1 \\
  \theta^2 \\
  \theta^3 \\
  \theta^4 \\
\end{array}%
\right)+
\left(%
\begin{array}{c}
  T^1_{24}~\theta^2 \wedge \theta^4 \\
  T^2_{34}~\theta^3 \wedge \theta^4  \\
  0 \\
  0 \\
\end{array}%
\right)\\
\end{equation}
The infinitesimal action on the torsion is
\begin{equation}\label{b4}
\left.
\begin{array}{cccc}
d~T^1_{24}\equiv(\alpha_1-\alpha_3-\alpha_8)T^1_{24}\\
d~T^2_{34}\equiv(\alpha_3-\alpha_6-\alpha_8)T^2_{34}\\
\end{array}
\right\}
~\textrm{mod}~(\theta^1,\theta^2,\theta^3,\theta^4)\\
\end{equation}
and a parametric calculation gives $T^1_{24}=-\frac{a_1}{a_3
a_8}\ne0$ and $T^2_{34}=-\frac{a_3}{a_6 a_8}\ne0$. We normalize
the torsion by setting
\begin{equation}\label{b5}
\begin{array}{cc}
T^1_{24}=-1,&T^2_{34}=-1.\\
\end{array}
\end{equation}
This leads to the principal components
\begin{equation}\label{b6}
\begin{array}{cc}
\alpha_6=2\alpha_3-\alpha_1, &\alpha_8=\alpha_1-\alpha_3.\\
\end{array}
\end{equation}
The normalizations force relations on the group $G$ in the form
\begin{equation}\label{b61}
\begin{array}{cc}
a_6=\frac{a_3^2}{a_1},&a_8=\frac{a_1}{a_3}.\\
\end{array}
\end{equation}
The {\it first-order} normalizations yield  an adapted coframe
with the {\it six-dimensional group} $G_1$
\begin{equation}\label{b62}
\left(%
\begin{array}{c}
  \theta^1 \\
  \theta^2 \\
  \theta^3 \\
  \theta^4 \\
\end{array}%
\right)=\left(%
\begin{array}{cccc}
  a_1 & 0 & 0 & 0 \\
  a_2 & a_3 & 0 & 0 \\
  a_4 & a_5 & \frac{a_3^2}{a_1} & 0 \\
  a_7 & 0 & 0 & \frac{a_1}{a_3} \\
\end{array}%
\right) \left(%
\begin{array}{c}
  \omega^1 \\
  \omega^2 \\
  \omega^3 \\
  \omega^4 \\
\end{array}%
\right).
\end{equation}
This leads to the structure equation
\begin{equation}\label{b7}
d\left(%
\begin{array}{c}
  \theta^1 \\
  \theta^2 \\
  \theta^3 \\
  \theta^4 \\
\end{array}%
\right)=\left(%
\begin{array}{cccc}
  \alpha_1 & 0 & 0 & 0 \\
  \alpha_2 & \alpha_3 & 0 & 0 \\
  \alpha_4 & \alpha_5 & 2\alpha_3-\alpha_1 & 0 \\
  \alpha_7 & 0 & 0 & \alpha_1-\alpha_3 \\
\end{array}%
\right)\wedge \left(%
\begin{array}{c}
  \theta^1 \\
  \theta^2 \\
  \theta^3 \\
  \theta^4 \\
\end{array}%
\right)+
\left(%
\begin{array}{c}
  -\theta^2 \wedge \theta^4 \\
 -\theta^3 \wedge \theta^4  \\
  T^3_{34}~\theta^3 \wedge \theta^4  \\
  0  \\
\end{array}%
\right)\\
\end{equation}
The infinitesimal action on the torsion is
\begin{equation}\label{b8}
\begin{array}{cccc}
d~T^3_{34}\equiv (\alpha_3-\alpha_1)T^3_{34}+3(\alpha_2-\alpha_5)\\
\end{array}
~\textrm{mod}~(\theta^1,\theta^2,\theta^3,\theta^4)\\
\end{equation}
and we can translate $T^3_{34}$ to zero:
\begin{equation}\label{b9}
\begin{array}{cc}
T^3_{34}=0.\\
\end{array}
\end{equation}
This leads to the principal components
\begin{equation}\label{b91}
\begin{array}{cc}
\alpha_5=\alpha_2.\\
\end{array}
\end{equation}
The normalizations force relations on the group $G_1$ in the form
\begin{equation}\label{b92}
\begin{array}{c}
a_5=\frac{a_2a_3}{a_1}-\frac{a_3^2}{3a_1} s_1,\\
\end{array}
\end{equation}
where $s_1=f_q.$

The {\it second-order} normalizations yield  an adapted coframe
with the {\it five-dimensional group} $G_2$
\begin{equation}\label{b93}
\left(%
\begin{array}{c}
  \theta^1 \\
  \theta^2 \\
  \theta^3 \\
  \theta^4 \\
\end{array}%
\right)=\left(%
\begin{array}{cccc}
  a_1 & 0 & 0 & 0 \\
  a_2 & a_3 & 0 & 0 \\
  a_4 & \frac{a_2a_3}{a_1}-\frac{a_3^2}{3a_1} s_1 & \frac{a_3^2}{a_1} & 0 \\
  a_7 & 0 & 0 & \frac{a_1}{a_3} \\
\end{array}%
\right) \left(%
\begin{array}{c}
  \omega^1 \\
  \omega^2 \\
  \omega^3 \\
  \omega^4 \\
\end{array}%
\right).
\end{equation}
This leads to the structure equation
\begin{equation}\label{b10}
d\left(%
\begin{array}{c}
  \theta^1 \\
  \theta^2 \\
  \theta^3 \\
  \theta^4 \\
\end{array}%
\right)=\left(%
\begin{array}{cccc}
  \alpha_1 & 0 & 0 & 0 \\
  \alpha_2 & \alpha_3 & 0 & 0 \\
  \alpha_4 & \alpha_2 & 2\alpha_3-\alpha_1 & 0 \\
  \alpha_7 & 0 & 0 & \alpha_1-\alpha_3 \\
\end{array}%
\right)\wedge \left(%
\begin{array}{c}
  \theta^1 \\
  \theta^2 \\
  \theta^3 \\
  \theta^4 \\
\end{array}%
\right)+
\left(%
\begin{array}{c}
  -\theta^2 \wedge \theta^4 \\
 -\theta^3 \wedge \theta^4  \\
  T^3_{24}~\theta^2 \wedge \theta^4  \\
  T^4_{24}~\theta^2 \wedge \theta^4  \\
\end{array}%
\right)\\
\end{equation}
The infinitesimal action on the torsion is
\begin{equation}\label{b11}
\left.
\begin{array}{ll}
d~T^3_{24}\equiv 2(\alpha_3-\alpha_1)T^3_{24}-2\alpha_4\\
d~T^4_{24}\equiv-\alpha_3~T^4_{24}-\alpha_7\\
\end{array}
\right\}
~\textrm{mod}~(\theta^1,\theta^2,\theta^3,\theta^4)\\
\end{equation}
and we can translate $T^3_{24}$ and $T^4_{24}$  to zero:
\begin{equation}\label{b12}
\begin{array}{cc}
T^3_{24}=0,&T^4_{24}=0.\\
\end{array}
\end{equation}
This leads to the principal components
\begin{equation}\label{b121}
\begin{array}{cc}
\alpha_4=0,&\alpha_7=0.\\
\end{array}
\end{equation}
The normalizations force relations on the group $G_2$ as
\begin{equation}\label{b122}
\begin{array}{cc}
a_4=\frac{a^2_2}{2a_1}-\frac{a^2_3}{18 a_1}s_2,& a_7=\frac{a_1}{6a_3}s_3,\\
\end{array}
\end{equation}
where $s_2= 2f_q ^2 +9\,f_p   - 3\,D_x f_q ,~s_3=f_{q,q}$.

The {\it third-order} normalizations yield  an adapted coframe
with the {\it three-dimensional group} $G_3$
\begin{equation}\label{b123}
\left(%
\begin{array}{c}
  \theta^1 \\
  \theta^2 \\
  \theta^3 \\
  \theta^4 \\
\end{array}%
\right)=\left(%
\begin{array}{cccc}
  a_1 & 0 & 0 & 0 \\
  a_2 & a_3 & 0 & 0 \\
  \frac{a^2_2}{2a_1}-\frac{a^2_3}{18 a_1}s_2 & \frac{a_2a_3}{a_1}-\frac{a_3^2}{3a_1} s_1 & \frac{a_3^2}{a_1} & 0 \\
  \frac{a_1}{6a_3}s_3 & 0 & 0 & \frac{a_1}{a_3} \\
\end{array}%
\right) \left(%
\begin{array}{c}
  \omega^1 \\
  \omega^2 \\
  \omega^3 \\
  \omega^4 \\
\end{array}%
\right).
\end{equation}
This gives rise to the structure equation
\begin{equation}\label{b13}
d\left(%
\begin{array}{c}
  \theta^1 \\
  \theta^2 \\
  \theta^3 \\
  \theta^4 \\
\end{array}%
\right)=\left(%
\begin{array}{cccc}
  \alpha_1 & 0 & 0 & 0 \\
  \alpha_2 & \alpha_3 & 0 & 0 \\
  0 & \alpha_2 & 2\alpha_3-\alpha_1 & 0 \\
  0 & 0 & 0 & \alpha_1-\alpha_3 \\
\end{array}%
\right)\wedge \left(%
\begin{array}{c}
  \theta^1 \\
  \theta^2 \\
  \theta^3 \\
  \theta^4 \\
\end{array}%
\right)+
\left(%
\begin{array}{c}
  -\theta^2 \wedge \theta^4 \\
 -\theta^3 \wedge \theta^4  \\
  T^3_{14}~\theta^1 \wedge \theta^4  \\
  T^4_{12}~\theta^1 \wedge \theta^2+T^4_{13}~\theta^1 \wedge \theta^3  \\
\end{array}%
\right)\\
\end{equation}
The infinitesimal action on the torsion is
\begin{equation}\label{b14}
\left.
\begin{array}{cc}
d~T^3_{14}\equiv &-3(\alpha_1-\alpha_3) T^3_{14}\\
d~T^4_{12}\equiv &-2 \alpha_3 T^4_{12}-\alpha_2 T^4_{13} \\
d~T^4_{13}\equiv &(\alpha_1-3\alpha_3) T^4_{13}\\
\end{array}
\right\}
~\textrm{mod}~(\theta^1,\theta^2,\theta^3,\theta^4)\\
\end{equation}
and here we have a bifurcation in the flowchart depending on
whether $T^3_{14}$, $T^4_{12}$ and $T^4_{13}$ are zero. A
parametric calculation gives
\begin{equation}\label{b15}
\begin{array}{ll}
T^3_{14}=&-\frac{a_3^3~I_3}{a_1^3},\\
T^4_{13}=& -\frac{a_1~I_1}{6~a_3^3}, \hfill \\
T^4_{12}=& -\frac{I_2}{36~a_3^2}~~\textrm{mod}~~T^4_{13}, \hfill \\
\end{array}
\end{equation}
where
\begin{equation}\label{b16}
\begin{array}{ll}
  I_1  &={s_3}_q= f_{q,q,q},\\
  I_2  &={s^2_3}+6{s_3}_p =f_{q,q} ^2  + 6\,f_{p,q,q},  \hfill \\
  I_3  &=J^3=\frac{1}{54} \left(2s_1s_2-3D_x s_2+54f_u\right)\\
       &= \frac{1}{54} \left( {4\,f_q ^3  + 18\,f_q \left( {f_p  - D_x f_q } \right) + 9\,D_x^2 f_q  - 27\,D_x f_p + 54\,f_u }\right).  \hfill \\
\end{array}
\end{equation}
\section*{Branch 1. $I_1=I_2=0$ and $I_3\ne0$.}
We choose this branch since the third-order ODEs with five point
symmetries have the canonical form $u''' = s~u'+u$, where $s$ is
constant. Clearly,  $I_1=I_2=0$ and $I_3=1$  for this canonical
form.

We normalize the torsion $T^3_{14}$ by setting
\begin{equation}\label{b17}
\begin{array}{cc}
T^3_{14}=-1.\\
\end{array}
\end{equation}
This leads to the principal components
\begin{equation}\label{b171}
\begin{array}{c}
\alpha_3=\alpha_1.\\
\end{array}
\end{equation}
The normalizations force relations on the group $G_3$ in the form
\begin{equation}\label{b172}
\begin{array}{c}
a_3=\frac{a_1}{J}.\\
\end{array}
\end{equation}
The {\it fourth-order} normalizations yield  an adapted coframe
with the {\it two-dimensional group} $G_4$
\begin{equation}\label{b173}
\left(%
\begin{array}{c}
  \theta^1 \\
  \theta^2 \\
  \theta^3 \\
  \theta^4 \\
\end{array}%
\right)=\left(%
\begin{array}{cccc}
  a_1 & 0 & 0 & 0 \\
  a_2 & \frac{a_1}{J} & 0 & 0 \\
  \frac{a^2_2}{2a_1}-\frac{a_1}{18J^2}s_2 & \frac{a_2}{J}-\frac{a_1}{3J^2} s_1 & \frac{a_1}{J^2} & 0 \\
  \frac{1}{6}Js_3 & 0 & 0 & J \\
\end{array}%
\right) \left(%
\begin{array}{c}
  \omega^1 \\
  \omega^2 \\
  \omega^3 \\
  \omega^4 \\
\end{array}%
\right).
\end{equation}
This leads to the structure equation
\begin{equation}\label{b18}
d\left(%
\begin{array}{c}
  \theta^1 \\
  \theta^2 \\
  \theta^3 \\
  \theta^4 \\
\end{array}%
\right)=\left(%
\begin{array}{cccc}
  \alpha_1 & 0 & 0 & 0 \\
  \alpha_2 & \alpha_1 & 0 & 0 \\
  0 & \alpha_2 & \alpha_1 & 0 \\
  0 & 0 & 0 & 0 \\
\end{array}%
\right)\wedge \left(%
\begin{array}{c}
  \theta^1 \\
  \theta^2 \\
  \theta^3 \\
  \theta^4 \\
\end{array}%
\right)+
\left(%
\begin{array}{c}
  -\theta^2 \wedge \theta^4 \\
 T^2_{23}~\theta^2 \wedge \theta^3+T^2_{24}~\theta^2 \wedge \theta^4-\theta^3 \wedge \theta^4  \\
 -\theta^1 \wedge \theta^4+T^3_{23}~\theta^2 \wedge \theta^3+2~T^2_{24}~\theta^3 \wedge \theta^4  \\
 T^4_{14}~\theta^1 \wedge \theta^4-\frac{1}{2}T^3_{23}~\theta^2 \wedge \theta^4+T^2_{23}~\theta^3 \wedge \theta^4 \\
\end{array}%
\right)\\
\end{equation}
The infinitesimal action on the torsion is
\begin{equation}\label{b19}
\left.
\begin{array}{ll}
d~T^2_{23}\equiv &-\alpha_1~T^2_{23}\\
d~T^2_{24}\equiv &-\alpha_2\\
d~T^3_{23}\equiv &-\alpha_1~T^3_{23}+2\alpha_2~T^2_{23}\\
d~T^4_{14}\equiv &-\alpha_1~T^4_{14}+\frac{1}{2}\alpha_2~T^3_{23}\\
\end{array}
\right\}
~\textrm{mod}~(\theta^1,\theta^2,\theta^3,\theta^4)\\
\end{equation}
and here we have a bifurcation in the flowchart depending on
whether $T^2_{23}$, $T^3_{23}$ and  $T^4_{14}$ are zero. A
parametric calculation provides
\begin{equation}\label{b20}
\begin{array}{ll}
T^2_{23}=&\frac{J I_4}{a_1},\\
T^2_{24}=&-\frac{a_2}{a_1}+\frac{1}{3}\frac{1}{J^2}s_4,\\
T^3_{23}=&\frac{I_5}{3a_1}  ~~\textrm{mod}~~T^2_{23},\hfill \\
T^4_{14}=&\frac{I_6}{a_1J} ~~\textrm{mod}~~(T^2_{23},T^3_{23}), \hfill \\
\end{array}
\end{equation}
where
\begin{equation}\label{b21}
\begin{array}{ll}
  I_4  =J_q,\\
  I_5  = f_{q,q} \,J - 6\,J_p, \\
  I_6  = J_u  - D_x J_p ,  \hfill \\
  s_4= 3 D_x J-J f_q.\\
\end{array}
\end{equation}
\section*{Branch 1.1. $I_4=I_5=I_6=0$.}
Similarly, we choose this branch as $I_4=I_5=I_6=0$   for the
canonical form $u''' = s~u'+u$, where $s$ is constant.

We can translate $T^2_{24}$   to zero:
\begin{equation}\label{b22}
\begin{array}{c}
T^2_{24}=0.\\
\end{array}
\end{equation}
This yields the principal components
\begin{equation}\label{b221}
\begin{array}{c}
\alpha_2=0.\\
\end{array}
\end{equation}
The normalizations force relations on the group $G_4$ in the form
\begin{equation}\label{b222}
\begin{array}{c}
a_2=\frac{1}{3}\frac{a_1}{J^2}s_4.\\
\end{array}
\end{equation}
The {\it fifth-order} normalizations give  an adapted coframe
with the {\it one-dimensional group} $G_5$
\begin{equation}\label{b223}
\left(%
\begin{array}{c}
  \theta^1 \\
  \theta^2 \\
  \theta^3 \\
  \theta^4 \\
\end{array}%
\right)=\left(%
\begin{array}{cccc}
  a_1 & 0 & 0 & 0 \\
  \frac{1}{3}\frac{a_1}{J^2}s_4 & \frac{a_1}{J} & 0 & 0 \\
  \frac{1}{18}\frac{a_1}{J^4}s_4^2-\frac{a_1}{18J^2}s_2 & \frac{1}{3}\frac{a_1}{J^3}s_4-\frac{a_1}{3J^2} s_1 & \frac{a_1}{J^2} & 0 \\
  \frac{1}{6}Js_3 & 0 & 0 & J \\
\end{array}%
\right)  \left(%
\begin{array}{c}
  \omega^1 \\
  \omega^2 \\
  \omega^3 \\
  \omega^4 \\
\end{array}%
\right).
\end{equation}
This results in the structure equation
\begin{equation}\label{b23}
d\left(%
\begin{array}{c}
  \theta^1 \\
  \theta^2 \\
  \theta^3 \\
  \theta^4 \\
\end{array}%
\right)=\left(%
\begin{array}{cccc}
  \alpha_1 & 0 & 0 & 0 \\
  0 & \alpha_1 & 0 & 0 \\
  0 & 0 & \alpha_1 & 0 \\
  0 & 0 & 0 & 0 \\
\end{array}%
\right)\wedge \left(%
\begin{array}{c}
  \theta^1 \\
  \theta^2 \\
  \theta^3 \\
  \theta^4 \\
\end{array}%
\right)+
\left(%
\begin{array}{c}
  -\theta^2 \wedge \theta^4 \\
 T^2_{14}~\theta^1 \wedge \theta^4-\theta^3 \wedge \theta^4  \\
 T^3_{12}~\theta^1 \wedge \theta^2-\theta^1 \wedge \theta^4+T^2_{14}~\theta^2 \wedge \theta^4  \\
 0 \\
\end{array}%
\right)\\
\end{equation}
The infinitesimal  action on the torsion is
\begin{equation}\label{b24}
\left.
\begin{array}{ll}
d~T^2_{14}\equiv &0\\
d~T^3_{12}\equiv &-\alpha_1 T^3_{12}\\
\end{array}
\right\}
~\textrm{mod}~(\theta^1,\theta^2,\theta^3,\theta^4)\\
\end{equation}
and in this case we have a bifurcation in the flowchart depending on the
value of  $T^2_{14}$ and whether $T^3_{12}$ is zero. A parametric
calculation provides
\begin{equation}\label{b25}
\begin{array}{ll}
T^2_{14}=-\frac{I_8}{2~J^4},\\
T^3_{12}=-\frac{I_7}{9a_1~J},\\
\end{array}
\end{equation}
where
\begin{equation}\label{b26}
\begin{array}{ll}
  I_7  = f_{q,q} \left(f_q ^2+ 9\,f_p    - 3\,D_x f_q  \right) - 9\,f_{p,p}  + 18\,f_{u,q}  - 6\,f_q f_{p,q},  \hfill \\
  I_8 = \frac{1}{3} \left( \left( {f_q ^2  + 3\,f_p  - 3\,D_x f_q \,} \right)\,J^2 \, + 6\,J\,D_x^2 J - 9\,\left( {D_x J} \right)^2\right) . \\
\end{array}
\end{equation}
It should be noted here that the relative invariant  $I_7=0$ and
the invariant $\frac{I_8}{J^4}=s$ for the canonical form $u''' =
s~u'+u$. Therefore, we choose the branch $I_7=0,
\frac{I_8}{J^4}=K$ where $K$ is constant.  Equivalently, we choose
the branch $I_7=K_q=K_p=K_u=K_x=0$.
\section*{Branch 1.1.1. $I_7=K_q=K_p=K_u=K_x=0$.}
In this branch, there is no more unabsorbable torsion left, so the
final remaining group variable $a_1$ cannot be normalized.
Moreover, $\alpha_1$ is now uniquely defined, so the problem is
determinant. This results in the following e-structure on the
five-dimensional prolonged space $M^{(1)}=M \times G_5$
\begin{equation}\label{b36}
\begin{array}{llll}
  \theta^1=a_1\omega^1,\\
  \theta^2=\frac{1}{3}\frac{a_1}{J^2}s_4\omega^1+\frac{a_1}{J}\omega^2,\\
  \theta^3=\frac{1}{18}(\frac{a_1}{J^4}s_4^2-\frac{a_1}{J^2}s_2)\omega^1+ \frac{1}{3}(\frac{a_1}{J^3}s_4-\frac{a_1}{J^2} s_1)\omega^2+\frac{a_1}{J^2}\omega^3,\\
  \theta^4=\frac{1}{6}Js_3\omega^1+J \omega^4,\\
  \alpha_1=\frac{da_1}{a_1}+\frac{1}{36 J}(4s_3 s_4-18 {s_4}_p+J s_1 s_3-6J {s_1}_p)\omega^1+\frac{1}{6}s_3 \omega^2-\frac{1}{3}\frac{s_4}{J}\omega^4.\\
\end{array}
\end{equation}
This results in the structure equations
\begin{equation}\label{b35}
\begin{array}{l}
  d\theta^1=-\theta^1 \wedge \alpha_1-\theta^2 \wedge \theta^4 \\
  d\theta^2=-\frac{K}{2}~\theta^1 \wedge \theta^4-\theta^2 \wedge \alpha_1-\theta^3 \wedge \theta^4  \\
  d\theta^3=-\theta^1 \wedge \theta^4-\frac{K}{2}~\theta^2 \wedge \theta^4-\theta^3 \wedge \alpha_1  \\
  d\theta^4=0 \\
  d\alpha_1=0 \\
\end{array}
\end{equation}
The invariant structure of the prolonged coframe are all constant.
We have produced an invariant coframe with rank zero on the
five-dimensional space coordinates $x,u,p, q, a_1$. Any such
differential equation admits a five-dimensional symmetry group of
point transformations.

Moreover, the symmetrical version of the Cartan formulation
$\textrm{mod}~(\omega^1,\omega^2,\omega^3)$ is
\begin{equation}\label{b37}
\begin{array}{l}
  J dx=\bar{J} d\bar{x} \\
  \frac{d a_1}{a_1}-\frac{1}{3}\frac{s_4}{J}dx=\frac{d \bar{a}_1}{\bar{a}_1}-\frac{1}{3}\frac{\bar{s}_4}{\bar{J}}d\bar{x}.\\
\end{array}
\end{equation}
Inserting the point transformation (\ref{ccc}) into (\ref{b37})
and using  $\bar{J}=1,~\bar{s}_4=0$ for
$\bar{f}=s\bar{u}'+\bar{u}$ and $\bar{a}=1$, results in
\begin{equation}\label{b39}
\begin{array}{l}
  J=D_x \phi, \\
  \frac{D_x a_1}{a_1}=\frac{1}{3}\left(\frac{3 D_x J-J f_q}{J}\right),\\
\end{array}
\end{equation}
where the auxiliary function $a_1(x,u,p)=\frac{\phi_x\psi_u   -
\phi_u \psi_x}{D_x \phi}$. This proves the following theorem.
\begin{theorem}
The necessary and sufficient conditions for equivalence of a
third-order equation $ u''' = f(x,u,u',u'')$ to the canonical form
$\bar{u}'''=s~\bar{u}'+\bar{u},~s=\,{\rm constant}$, with five point
symmetries via {\it point transformation} (\ref{ccc}) are the
identically vanishing of the relative invariants
\begin{equation}\label{b40}
\begin{array}{l}
  I_1  = f_{q,q,q}  \hfill \\
  I_2  = f_{q,q} ^2  + 6\,f_{p,q,q}  \hfill \\
  I_4  = J_q \,\,\,\,\,\,\,\,\,\,\,\,\,\,\,\,\,\,\,\,\,\,\,\,\,\,\,\,\,\,\,\,\,\,\,\,\,\,\,\,\,\,\,\,\,\,\,\,\,\,\,\,\,\,\,\,\,\,\,\,\,\,\,\,\,\,\,\,\,\,\,\,\,\,\,\,\,\,\,\,\,\,\,\,\,\,\,\,\,\,\,\,\,\,\,\,\,\,\,\,\, \hfill \\
  I_5  = f_{q,q} \,J - 6\,J_p \,\,\,\,\,\,\,\,\,\,\,\,\,\,\,\,\,\,\,\,\,\,\,\,\,\,\,\,\,\,\, \hfill \\
  I_6  = J_u  - D_x J_p  \hfill \\
  I_7  = f_{q,q} \left( {9\,f_p  + f_q ^2  - 3\,D_x f_q } \right) - 9\,f_{p,p}  + 18\,f_{u,q}  - 6\,f_q f_{p,q}  \hfill \\
  I_9  =K_q\\
  I_{10}=K_p\\
  I_{11}=K_u\\
  I_{12}=K_x,\\
\end{array}
\end{equation}
where
\begin{equation}\label{b41}
\begin{array}{l}
I_3  =J^3= \frac{1}{54} \left( {4\,f_q ^3  + 18\,f_q \left( {f_p  - D_x f_q } \right) + 9\,D_x^2 f_q    - 27\,D_x f_p+ 54\,f_u  }\right)  \ne 0,\\
I_8 = \frac{1}{3} \left( \left( {f_q ^2  + 3\,f_p  - 3\,D_x f_q \,} \right)\,J^2 \, + 6\,J\,D_x^2 J - 9\,\left( {D_x J} \right)^2\right), \\
K=\frac{I_8}{J^4}.\\
\end{array}
\end{equation}
Given that the the system of relative invariants (\ref{b40}) is
zero, the linearizing point transformation (\ref{ccc}) is defined by
\begin{equation}\label{b42}
\begin{array}{l}
 J=D_x \phi, \\
  \frac{D_x a_1}{a_1}=\frac{1}{3}\left(\frac{3 D_x J-J f_q}{J}\right),\\
\end{array}
\end{equation}
where the auxiliary function $a_1(x,u,p)=\frac{1}{J}(\phi_x\psi_u
- \phi_u \psi_x)$.

Finally, the constant $s$ of the resulting canonical form is given
by the equation $s=K$.
\end{theorem}
\section{Illustration of the theorem}
\begin{example}\rm \cite{ibr}
Consider the nonlinear ODE
\begin{equation}\label{n1}
\begin{array}{c}
u'''=\frac{3{u''}^2}{u'}+x{u'}^4.\\
\end{array}
\end{equation}
The function
\begin{equation}\label{n2}
\begin{array}{ll}
f(x,u,p,q)=\frac{3{q}^2}{p}+x{p}^4\\
\end{array}
\end{equation} satisfies the constraints $I_1=I_2=I_4=I_5=I_6=I_7=I_9=I_{10}=I_{11}=I_{12}=0$ while $I_3\ne0$;
consequently, this equation admits the five-dimensional point
symmetry group. Moreover, since $I_8=0$ and $K=\frac{I_8}{J^4}=0$,
then it is equivalent to the canonical form $\bar{u}'''=\bar{u}.$

Since $J=-p$, then the linearizing transformation (\ref{ccc}) can
be obtained by (\ref{b42}) as
\begin{equation}\label{n3}
\begin{array}{lll}
  \phi_x=0,&\phi_u=-1 \\
  \frac{1}{a_1}\left(\frac{\partial a_1}{\partial x}+\frac{\partial a_1}{\partial u}p\right)=0,&\frac{1}{a_1}\frac{\partial a_1}{\partial p}=-\frac{1}{p}.\\
\end{array}
\end{equation}
A solution of the system (\ref{n3}) is $\phi=-u$ and
$a_1=\frac{1}{p}$. Since, the auxiliary function
$a_1=\frac{\psi_x}{p}$, then $\psi=x$. Therefore, the canonical
form $\bar{u}'''=\bar{u}$  can be obtained for the ODE (\ref{n1}) via
the point transformation $$\bar{x}=-u, \bar{u}=x.$$
\end{example}
\begin{example}\rm
We now focus on the linear ODE with variable coefficients
\begin{equation}\label{n4}
\begin{array}{c}
u'''=-\frac{3}{x}u''+(8x^2+\frac{3}{x^2})u'+8x(x^2+2)u.\\
\end{array}
\end{equation}
The function
\begin{equation}\label{n5}
\begin{array}{ll}
f(x,u,p,q)=-\frac{3}{x}q+(8x^2+\frac{3}{x^2})p+8x(x^2+2)u\\
\end{array}
\end{equation} satisfies the constraints $I_1=I_2=I_4=I_5=I_6=I_7=I_9=I_{10}=I_{11}=I_{12}=0$ whereas $I_3\ne0$;
consequently, this equation admits the five-dimensional point
symmetry group. Moreover, since $J=2x,~I_8=32x^4$ and
$K=\frac{I_8}{J^4}=2$, then it is equivalent to the canonical form
$\bar{u}'''=2\bar{u}'+\bar{u}$.

The linearizing transformation can be obtained via (\ref{b42}) by
solving the following system
\begin{equation}\label{n5}
\begin{array}{lll}
  \phi_x=2x,&\phi_u=0 \\
  \frac{1}{a_1}\left(\frac{\partial a_1}{\partial x}+\frac{\partial a_1}{\partial u}p\right)=\frac{2}{x},&\frac{\partial a_1}{\partial p}=0.\\
\end{array}
\end{equation}
A solution of the system (\ref{n5}) is $\phi=x^2$ and $a_1=x^2$.
Since, the auxiliary function $a_1=\psi_u$, then $\psi=x^2 u$.
Therefore, the canonical form $\bar{u}'''=2\bar{u}'+\bar{u}$  can
be achieved for the ODE (\ref{n4}) via the transformation
$$\bar{x}=x^2, \bar{u}=x^2u.$$
\end{example}
\begin{example}\rm
Consider now the linear ODE with variable coefficients in the
Laguerre-Forsyth canonical form
\begin{equation}\label{n6}
\begin{array}{c}
u'''=\frac{1}{x^6}u.\\
\end{array}
\end{equation}
The function
\begin{equation}\label{n7}
\begin{array}{ll}
f(x,u,p,q)=\frac{u}{x^6}\\
\end{array}
\end{equation} satisfies the constraints $I_1=I_2=I_4=I_5=I_6=I_7=I_9=I_{10}=I_{11}=I_{12}=0$ whereas $I_3\ne0$;
consequently, this equation admits the five-dimensional point
symmetry group. Moreover, since $J=\frac{1}{x^2},~I_8=0$ and
$K=\frac{I_8}{J^4}=0$, then it is equivalent to the canonical form
$\bar{u}'''=\bar{u}$.

The linearizing transformation can be obtained via (\ref{b42}) by
solving the following system
\begin{equation}\label{n8}
\begin{array}{lll}
  \phi_x=\frac{1}{x^2},&\phi_u=0 \\
  \frac{1}{a_1}\left(\frac{\partial a_1}{\partial x}+\frac{\partial a_1}{\partial u}p\right)=-\frac{2}{x},&\frac{\partial a_1}{\partial p}=0.\\
\end{array}
\end{equation}
A solution of the system (\ref{n8}) is $\phi=-\frac{1}{x}$ and
$a_1=\frac{1}{x^2}$. Since, the auxiliary function $a_1=\psi_u$,
then $\psi=\frac{u}{x^2}$. Therefore, the canonical form
$\bar{u}'''=\bar{u}$ can be achieved for the ODE (\ref{n6}) via
the transformation
$$\bar{x}=-\frac{1}{x}, \bar{u}=\frac{u}{x^2}.$$
\end{example}
\begin{example}\rm
Consider now the nonlinear ODE
\begin{equation}\label{n9}
\begin{array}{c}
u'''=\frac{3}{2}\frac{{u''}^2}{u'}.\\
\end{array}
\end{equation}
The function
\begin{equation}\label{n10}
\begin{array}{ll}
f(x,u,p,q)=\frac{3}{2}\frac{{q}^2}{p}\\
\end{array}
\end{equation} gives $I_3=0$; consequently, this equation can not admit five-dimensional point
symmetry group. This agrees with the fact that the ODE (\ref{n9})
has six point symmetries \cite{IbrFazal1996}.
\end{example}
\section{Conclusion}
We have invoked the Cartan equivalence method to effectively and compactly solve the linearization
problem for a scalar third-order ODE to enable its reduction to a linear third-order ODE with five point symmetries.
Moreover, we were able to obtain the point transformation that does the reduction to this canonical form.
In previous work as in \cite{ibr}, inter alia, the authors utilize the direct method, to find reduction to the
Laguerre-Forsyth canonical form which can have four, five or seven point symmetries so there isn't a unique form. Notwithstanding, the five symmetry
case in the Laguerre-Forsyth canonical form as we have pointed out can have variable coefficient while in the canonical we
have utilized the coefficient is constant and the form simpler.

It is important to also mention herein that in the Cartan approach used, we have for the first time provided
how one can deduce the point transformations that does the reduction to the linear canonical form. This was previously done for the
direct method in \cite{ibr}.

Amongst the basic approaches to the linearization problem via point transformation are the two prominent methods of Lie
and Cartan. The first has been nicely utilized in \cite{ibr} and here we have effectively invoked the Cartan approach
for the five symmetry case. It remains to pursue the four symmetry case. The reader is referred to \cite{Dweik2016}
for the maximal symmetry case wherein compact criteria is provided.

We also have amply demonstrated our results by means of examples.

\subsection*{Acknowledgments}
Ahmad Y. Al-Dweik would like to thank the King Fahd University of
Petroleum and Minerals for its support and excellent research
facilities. FMM thanks the NRF of South Africa for research support.

\end{document}